\documentclass{amsart}
\usepackage{amssymb}
\usepackage{amsfonts}
\usepackage{latexsym}

\newtheorem{theorem}{Theorem}[section]
\newtheorem{lemma}[theorem]{Lemma}
\newtheorem{proposition}[theorem]{Proposition}
\newtheorem{corollary}[theorem]{Corollary}
\theoremstyle{definition}
\newtheorem{definition}[theorem]{Definition}
\newtheorem{example}[theorem]{Example}

\newtheorem{remark}[theorem]{Remark}

\newcommand{\id}{\text{id}}

\newcommand{\Hom}{\text{Hom}}

\newcommand{\eps}{\varepsilon}

\newcommand{\I}{\mathcal{I}}

\newcommand{\actl}{\rightharpoonup}
\newcommand{\actr}{\leftharpoonup}
\newcommand{\lact}{\triangleright}
\newcommand{\ract}{\triangleleft}
\newcommand{\la}{\langle\,}
\newcommand{\ra}{\,\rangle}

\newcommand{\0}{_{(0)}}
\newcommand{\1}{_{(1)}}
\newcommand{\2}{_{(2)}}
\newcommand{\3}{_{(3)}}
\newcommand{\4}{_{(4)}}
\renewcommand{\I}{^{(1)}}
\newcommand{\II}{^{(2)}}

\newcommand{\ot}{\otimes}

\begin{document}

\title{On twisting of finite-dimensional Hopf algebras}

\author{Eli Aljadeff}
\address{Department of Mathematics, Technion-Israel Institute of
Technology, Haifa 32000, Israel}
\email{aljadeff@math.technion.ac.il}

\author{Pavel Etingof}
\address{Department of Mathematics, Massachusetts Institute of Technology,
Cambridge, MA 02139, USA}
\email{etingof@math.mit.edu}

\author{Shlomo Gelaki}
\address{Department of Mathematics, Technion-Israel Institute of
Technology, Haifa 32000, Israel}
\email{gelaki@math.technion.ac.il}

\author{Dmitri Nikshych}
\address{Department of Mathematics and Statistics,
University of New Hampshire,  Durham, NH 03824, USA}
\email{nikshych@math.unh.edu}

\begin{abstract}
In this paper we study
the properties of Drinfeld's twisting for finite-dimensional Hopf
algebras. We determine how the integral of
the dual to a unimodular Hopf algebra $H$ changes under twisting
of $H$. We show that the classes of cosemisimple unimodular, cosemisimple
involutive, cosemisimple quasitriangular finite-dimensional Hopf algebras
are stable under twisting. We also prove the cosemisimplicity of a coalgebra
obtained by twisting of a cosemisimple unimodular Hopf algebra
by two different twists on two sides (such twists are closely
related to biGalois extensions \cite{Sc}), and describe the
representation theory of its dual. Next, we define the notion of a
non-degenerate twist for a Hopf algebra $H$, and set up a bijection
between such twists for $H$ and $H^*$.
This bijection is based on Miyashita-Ulbrich actions
of Hopf algebras on simple algebras \cite{DT2, U}.
It generalizes to the non-commutative case the procedure of
inverting a non-degenerate skew-symmetric bilinear form on a vector space.
Finally, we apply these results to classification of
twists in group algebras and of cosemisimple triangular finite-dimensional
Hopf algebras in positive characteristic, generalizing the
previously known classification in characteristic zero.
\end{abstract}
\maketitle

%%%%%%%%%%%%%%%%%%%%%%%%%%%%%%%%%%%%%%%%%%%%%%%%%%%%%%%%%%%%%%%%%%%%%%%%%%%%%%%%%%%%%%%%%%%%%%%%%%%%%%%%%%%%%%%%%%%%%%%%%%%%%%%%%%%%%%%%%%%%%%%%%%%%%%%%

\begin{section}
{Introduction} In this paper we study some properties of
Drinfeld's twisting operation for finite-dimensional Hopf
algebras. In particular, for a finite-dimensional unimodular Hopf
algebra $H$, we determine how the integral of the dual Hopf
algebra $H^*$ changes under twisting. This allows us to show that
the classes of cosemisimple unimodular, cosemisimple involutive,
cosemisimple quasitriangular finite-dimensional Hopf algebras are
invariant under twisting for any characteristic of the ground
field $k$.

We also consider coalgebras obtained from $H$ by twisting it using
two different twists,  on the left and on the right (this reduces
to Hopf algebra twisting if these two twists  are the same). We
show that such a coalgebra is always cosemisimple if $H$ is
cosemisimple and unimodular. In particular, this applies to the
situation when one  of the two twists equals to $1$ (twisting on
one side). In this case, the unimodularity assumption can be
dropped, and this result is a known theorem of Blattner and
Montgomery. We also extend the result of \cite{EG2} by describing the algebra
structure of the dual of the coalgebra obtained by twisting
of  $\mathbb{C}[G]$, the group algebra of a finite group $G$,
by means of two different twists.
%These results generalize to the case of two pseudo-twists, if they define
%the same associator.

Twisting $H$ on one side defines not simply a coalgebra but
actually an $H$-module coalgebra $C$ which is the regular
$H$-module, such that the $H$-module algebra extension $k\subset
C^*$ is Galois. It is easy to show that the converse is also true:
any $H$-module coalgebra with such properties comes from twisting
$H$ on one side by a twist, which is unique up to a gauge
transformation. Thus, our results  imply that if $H$ is
cosemisimple, then any $H$-module coalgebra with
said properties is cosemisimple.

Next, we define a non-degenerate twist for $H$ as a twist for
which the corresponding  $H$-coalgebra is simple (this generalizes
the notion of a minimal twist for group algebras). The above
results, combined with Masuoka's generalization of the
Skolem-Noether theorem \cite{Mas} and the technique of
Miyashita-Ulbrich actions \cite{DT2, U}, allow us to show that
there exists a natural bijection between the set $NT(H)$ of gauge
classes of non-degenerate twists for $H$, and the set $NT(H^*)$;
in particular, these sets (which are finite in the semisimple
cosemisimple case by a theorem of H.-J.~Schneider \cite{Sch}) have
the same number of elements. For example, if $H=\mathbb{C}[V]$,
where $V$ is a vector space over a finite field of odd
characteristic, then $NT(H)$ is the set of symplectic structures
on $V^*$, $NT(H^*)$ is the set of symplectic structures on $V$,
and our bijection is the usual inversion. More generally, if $G$
is a finite group and $H=\mathbb{C}[G]$ then this bijection is the
Movshev bijection \cite{Mov} between gauge classes of minimal
twists for $\mathbb{C}[G]$ and cohomology classes of
non-degenerate $2$-cocycles on $G$ with values in
$\mathbb{C}^\times$.

Finally, we apply the results of invariance of cosemisimplicity
under twisting to the problem of classification of twists in group
algebras, and of finite-dimensional cosemisimple triangular Hopf
algebras in any characteristic (over algebraically closed fields),
which generalizes the previous results of Movshev \cite{Mov} and
the second and third authors \cite{EG1} in characteristic zero.
The answer turns out to be very simple: the classification is the
same as in characteristic  zero, except that  the subgroup $K$
with an irreducible projective representation of dimension
$|K|^{1/2}$ has to be of order coprime to the characteristic. In
particular, this gives a Hopf algebraic proof of the known fact that a
group of central type in characteristic $p$ has order coprime to
$p$.

{\bf Acknowledgments.} P.E.\ and D.N.\ were
partially supported by the NSF grant DMS-9988796.
P.E.\ and D.N.\ partially conducted  their research for the
Clay Mathematics Institute as a Clay Mathematics Institute Prize Fellow
and Liftoff Mathematician, respectively.
D.N.\ thanks the Mathematics Department of MIT for  the warm hospitality
during his visit. S.G.\ research was supported by the VPR - Fund at the
Technion and by the fund for the promotion of research at the Technion.

We are most grateful to the referee, for the explanation how to
use the Miyashita-Ulbrich actions and the results proved in the
papers \cite{DT1, DT2, Sc, U}, to strengthen several statements in
Section 5; in particular to drop the semisimplicity assumption
from Theorem~\ref{bijection for twists} and to show that
$D_{H^*}\circ D_H =\id$.
\end{section}

%%%%%%%%%%%%%%%%%%%%%%%%%%%%%%%%%%%%%%%%%%%%%%%%%%%%%%%%%%%%%%%%%%%%%%%%%%%%%

\begin{section}
{Preliminaries} Throughout this paper $k$ denotes an
algebraically closed field of
an arbitrary characteristic. We use Sweedler's notation for the
comultiplication in a coalgebra : $\Delta(c) = c\1 \otimes c\2,$
where summation is understood.

Let $H$ be a Hopf algebra (see \cite{M,Sw} for the definition)
over $k$ with comultiplication $\Delta$, counit $\eps$, and
antipode $S$. The following notion of a twisting deformation of
$H$ is due to V.~Drinfeld \cite{D}.

\begin{definition}
\label{definition of twist}
A {\em twist} for $H$ is an invertible element $J \in H\otimes H$
that satisfies
\begin{equation}
\label{twist eqn}
(\Delta\otimes \id)(J) (J \otimes 1) = (\id \otimes \Delta)(J) (1\otimes J).
\end{equation}
\end{definition}

\begin{remark}
\label{eps non-essential}
Applying $\eps$ to \eqref{twist eqn} one sees that
$c= (\eps\otimes\id)(J) =(\id\otimes\eps)(J)$ is a non-zero scalar
for any twist $J$, cf.\ \cite[p.~813]{DT1}.
In particular, one can always replace $J$ by $c^{-1}J$ to normalize
it in such a way that
\begin{equation}
\label{eps twist} (\eps\otimes\id)(J)=(\id\otimes\eps)(J)=1.
\end{equation}
We will always assume that $J$ is normalized in this way.
\end{remark}

We use a convenient shorthand notation for $J$ and $J^{-1}$,
writing $J = J\I \otimes J\II$ and $J^{-1} = J^{-(1)} \otimes
J^{-(2)}$, where of course a summation is understood. We also
write $J_{21}$ for $J^{(2)} \otimes J^{(1)}$ etc.

\begin{remark}
Let  $x\in H$ be an invertible element such that $\eps(x)= 1$. If
$J$ is a twist for $H$ then so is $J^x := \Delta(x) J (x^{-1}
\otimes x^{-1})$. The twists $J$ and $J^x$ are said to be {\em
gauge equivalent}.
\end{remark}

Given a twist $J$ for $H$  one can define a new Hopf algebra $H^J$
with the same algebra structure and counit as $H$, for which the
comultiplication and antipode are given by
\begin{eqnarray}
\Delta^J(h)&=& J^{-1} \Delta(h) J, \\ S^J(h) &=& Q_J^{-1} S(h)
Q_J, \qquad \mbox{ for all } \quad h\in H,
\end{eqnarray}
where $Q_J = S(J\I)J\II$ is an invertible element of $H$ with the inverse
$Q_J^{-1} = J^{-(1)} S(J^{-(2)})$.

The element $Q_J$ satisfies the following useful identity (cf.\
\cite[(2.17)]{Maj}) :
\begin{equation}
\label{Delta Q} \Delta(Q_J) = (S\otimes S)(J^{-1}_{21}) (Q_J
\otimes Q_J)  J^{-1}.
\end{equation}
In particular, we have,
\begin{equation}
\label{Q-1S(Q)} \Delta(Q_J^{-1}S(Q_J)) = J (Q_J^{-1}S(Q_J) \otimes
Q_J^{-1}S(Q_J)) (S^2\otimes S^2)(J^{-1}).
\end{equation}

Twists for Hopf algebras of finite groups were studied in \cite{Mov}
and \cite{EG1} in the case when
$\mbox{char}(k)$ is prime to the order
of the group. Such twists were
classified up to a gauge equivalence and were used  in \cite{EG1}
to classify all semisimple cosemisimple triangular Hopf algebras
over $k$.
We refer the reader to the survey \cite{G2} for a review of the recent
developments in the classification of triangular Hopf algebras.

In subsequent sections we will use the following lemma which is a direct
consequence of the twist equation~\eqref{twist eqn}.

\begin{lemma}
\label{twist identities}
Let $J$ be a twist for $H$. Then
\begin{eqnarray*}
S(J\I){J\II}\1 \otimes {J\II}\2 &=& (Q_J \otimes 1) J^{-1},\\
S({J\I}\1) \otimes S({J\I}\2) J\II &=& (S\otimes S)(J^{-1}) (1\otimes Q_J)\\
J^{-(1)}\1 \otimes J^{-(1)}\2 S(J^{-(2)}) &=& J (1\otimes Q_J^{-1}), \\
J^{-(1)} S(J^{-(2)}\1) \otimes S(J^{-(2)}\2) &=&
( Q_J^{-1}\otimes 1)(S\otimes S)(J).
\end{eqnarray*}
\end{lemma}

\end{section}

%%%%%%%%%%%%%%%%%%%%%%%%%%%%%%%%%%%%%%%%%%%%%%%%%%%%%%%%%%%%%%%%%%%%%%%%%%%%%

\begin{section}
{Twisting of unimodular Hopf algebras}

Let $H$ be a unimodular Hopf algebra and let $J$ be a twist for $H$.
Let $\lambda\in H^*$ (respectively, $\rho\in H^*$) be a non-zero
left (respectively, right) integral on $H$, i.e.,
\begin{equation*}
h\1 \la \lambda,\, h\2 \ra = \la \lambda,\, h\ra 1, \qquad
\la \rho,\, h\1 \ra h\2 = \la \rho,\, h \ra 1, \quad \mbox{ for all } h\in H.
\end{equation*}

\begin{remark}
\label{nakayama} Note that in the case of a finite-dimensional $H$,
by \cite[Theorem 3]{R}, for all $g, h
\in H$ we have
\begin{equation*}
\la \lambda,\, gh \ra = \la \lambda, h S^2(g) \ra, \qquad \la
\rho,\, gh \ra = \la \rho, S^2(h)g \ra.
\end{equation*}
\end{remark}

\begin{remark}
Let $H$ be an arbitrary  Hopf algebra and let $\lambda$ and $\rho$
be as above. We will use the following invariance properties of
integrals on a Hopf algebra, which are straightforward to check :
\begin{eqnarray}
\label{left invariance} g\1 \la \lambda,\,h g\2 \ra &=& S(h\1) \la
\lambda,\,h\2 g \ra,\\ \label{right invariance} \la \rho,\,g\1 h
\ra g\2 &=& \la \rho,\, g h\1 \ra S(h\2),
\end{eqnarray}
for all $g,h \in H$.
\end{remark}

\begin{remark}
Let $H$ be an arbitrary finite-dimensional Hopf algebra and let
$\lambda$ be as above. Recall that $H$ is a Frobenius algebra with
non-degenerate bilinear form $H\ot H\rightarrow k$ given by $h\ot
g\mapsto \la \lambda,\,hg \ra$ (cf. \cite{M}, \cite{Sw}).
\end{remark}

Our first goal is to describe the integrals on the twisted Hopf
algebra $H^J$ in the case when $H$ is finite-dimensional.

Recall that $H$ acts on $H^*$ on the left via $h \actl \phi =
\phi\1 \la \phi\2, \, h\ra$ and on the right via $\phi \actr h =
\la \phi\1, \, h\ra \phi\2$.

Let us denote $u_{_J} := Q_J^{-1}S(Q_J)$.

\begin{theorem}
\label{integrals in H_J} Let $H$ be a finite-dimensional
unimodular Hopf algebra and let $J$ be a twist for $H$. Let
$\lambda$ and $\rho$ be, respectively, non-zero left and right
integrals on $H$. Then the  elements $\lambda_J := u_{_J}\actl
\lambda$ and $\rho_J := \rho \actr u_{_J}^{-1}$ are, respectively,
non-zero left and right integrals on $H^J$.
\end{theorem}
\begin{proof}
We have
\begin{equation}
\label{WJ} \Delta(u_{_J}) = J (u_{_J} \otimes u_{_J})  (S^2\otimes
S^2)(J^{-1}).
\end{equation}
by equation~\eqref{Q-1S(Q)}. We need to check that
$S(h_{\tilde{(1)}}) \la \lambda_J,\,h_{\tilde{(2)}} \ra=\la
\lambda_J,\,h\ra1$ for all $h\in H^J$, where $\Delta^J(h) =
h_{\tilde{(1)}} \otimes h_{\tilde{(2)}}$. To this end, we
compute, using Remark~\ref{nakayama} and equation  \eqref{left
invariance} :
\begin{eqnarray*}
S(h_{\tilde{(1)}}) \la \lambda_J,\,h_{\tilde{(2)}} \ra &=&
S(J^{-(1)} h\1 J\I) \la \lambda,\, J^{-(2)} h\2 J\II u_{_J}\ra \\
&=& S(J\I) S(h\1) S(J^{-(1)}) \la \lambda,\, h\2 J\II u_{_J}
S^2(J^{-(2)})\ra \\ &=& S(J\I) {J\II}\1 {u_{_J}}\1
S^2({J^{-(2)}}\1) S(J^{-(1)})\\ & & \times \, \la \lambda,\, h
{J\II}\2 {u_{_J}}\2 S^2({J^{-(2)}}\2) \ra.
\end{eqnarray*}
Since $\lambda$ is non-degenerate, the identity in question is
equivalent to
\begin{equation*}
1\otimes u_{_J} = (S(J\I) {J\II}\1 \otimes {J\II}\2)\Delta(u_{_J})
(S^2({J^{-(2)}}\1) S(J^{-(1)}) \otimes S^2({J^{-(2)}}\2)),
\end{equation*}
which reduces to equation~\eqref{WJ} by  Lemma~\ref{twist identities}.

The proof of the statement regarding $\rho$ is completely similar.
\end{proof}

\begin{remark} Suppose that $H$ is cocommutative.
Then it is easy to show that $u_{_J}$ coincides with the Drinfeld
element $u$ of the triangular Hopf algebra $(H^J, J_{21}^{-1}J)$.
This motivated the notation $u_{_J}$.
\end{remark}

The next Corollary shows that the class of cosemisimple unimodular Hopf
algebras over $k$ is closed under twisting. This result is true regardless of
the characteristic of $k$ (in characteristic $0$ it easily follows
from the result of Larson and Radford \cite{LR}).

\begin{corollary}
\label{twisted Hopf alg} If $H$ is a cosemisimple  unimodular Hopf
algebra, then so is its twisting deformation $H^J$.
\end{corollary}
\begin{proof}
Clearly, $H^J$ is unimodular, since twisting preserves integrals
in $H$. Cosemisimplicity of $H^J$ is equivalent to $\la
\lambda_J,\, 1\ra \neq 0$ by Maschke's theorem. By
Theorem~\ref{integrals in H_J}, we have
\begin{eqnarray*}
\la \lambda_J,\, 1\ra &=& \la \lambda,\,Q_J^{-1}S(Q_J) \ra\\ &=&
\la \lambda,\, J^{-(1)} S(J^{-(2)}) S(J^{(2)}) S^2(J^{(1)}) \ra\\
&=& \la \lambda,\, J^{(1)} J^{-(1)} S(J^{-(2)}) S(J^{(2)}) \ra\\
&=& \la \lambda,\,1 \ra \ne 0
\end{eqnarray*}
where we used Remark~\ref{nakayama}.
\end{proof}

\begin{remark}
In \cite[Theorem 1.3.6]{G1} it was shown that a cosemisimple
quasitriangular Hopf algebra is automatically unimodular. Since
the property of being quasitriangular is preserved under a
twisting deformation, it follows from Corollary~\ref{twisted Hopf
alg} that the class of cosemisimple quasitriangular Hopf algebras
over a field of an arbitrary characteristic is closed under
twisting.
\end{remark}

\begin{remark}
\label{involutive} Suppose that $H$ is a finite-dimensional
cosemisimple involutive Hopf algebra, i.e., $S^2 =\id$.  Then
$u_{_J} = Q_J^{-1}S(Q_J)$ is a group-like element in $H^J$ by
Equation~\eqref{Q-1S(Q)}. Thus, we have that $ \lambda(1) =
\lambda_J(u_{_J}^{-1}) \neq 0$ (where $\lambda_J$ was defined in
Theorem \ref{integrals in H_J}) if and only if $u_{_J} =1$, i.e.,
if and only if $S(Q_J) = Q_J$. This means that $(S^J)^2 =\id$,
i.e., the class of finite-dimensional cosemisimple involutive Hopf
algebras over a field of an arbitrary characteristic is closed
under twisting.
\end{remark}

\begin{remark} We expect that the class of finite-dimensional cosemisimple
Hopf algebras is closed under twisting. This would follow from (a
weak form of) a conjecture of Kaplansky, saying that a
cosemisimple finite-dimensional Hopf algebra is unimodular.
\end{remark}

Below we prove a somewhat more general statement regarding the cosemisimplicity
of the coalgebra obtained by deforming $H$ by means of two different twists,
see Theorem~\ref{HLJ is cosemisimple} below. This result will be used in later
sections.

Recall \cite{Pi} that an algebra $A$ is said to be {\em separable} if
there exists an element $e = e\I \otimes e\II \in A\otimes A$,
where a summation is understood, such that
\begin{equation}
a e\I \otimes e\II = e\I \otimes e\II a \quad \mbox{ for all }
a\in A \quad \mbox{ and } \quad e\I e\II = 1.
\end{equation}
Such an element $e$ is called a {\em separability element} for
$A$. A separable algebra over a perfect field is always
finite-dimensional and semisimple \cite[10.2]{Pi}. For Hopf
algebras the notions of separability and  semisimplicity coincide;
moreover, a semisimple Hopf algebra is automatically
finite-dimensional.

Passing from algebras to coalgebras, one can define a notion
dual to separability as follows, cf.\ \cite{L}.

\begin{definition}
\label{coseparable}
A coalgebra $C$  is said to be {\em coseparable} if there exists
a bilinear form $\psi : C \otimes C \to k$ such that
\begin{equation}
\label{cosep pairing} c\1 \psi(c\2 \otimes d) = \psi(c \otimes
d\1)d\2  \quad\mbox{ for all } c,d\in C \quad \mbox{ and } \quad
\psi(\Delta(c)) = \eps(c).
\end{equation}
Such a form $\psi$ is called a {\em coseparability pairing} for $C$.
\end{definition}

\begin{remark}
\label{trivialities on cosep} For a coalgebra $C$ over an
algebraically closed field the notions of coseparability and
cosemisimplicity are equivalent.
\end{remark}

Let $J,L$ be twists for $H$, then
\begin{equation}
\label{2twist}
\tilde{\Delta}(h) := L^{-1} \Delta(h) J, \qquad h\in H
\end{equation}
defines a new coassociative comultiplication in $H$ for which
$\eps$ is still a counit. We will denote this twisted coalgebra by
$H^{(L,J)}$ and let $H^{(J)} = H^{(1,J)}$.

\begin{remark}
\label{closely related}
The two-sided twist \eqref{2twist}  is closely related to biGalois extensions
\cite{Sc}. The coalgebra $H^{(1, J)}$ (resp.\ $H^{(L, 1)}$, resp.\
$H^{(L, J)}$) is called $(H, H^J)$-biGalois (resp.\ $(H^L, H)$-biGalois,
resp.\ $(H^L, H^J)$-biGalois) coalgebra. One sees that $H^{(L,J)}
=H^{(L, 1)}\otimes_H H^{(1,J)}$.
\end{remark}

\begin{theorem}
\label{HLJ is cosemisimple}
Let $H$ be a finite-dimensional
cosemisimple  unimodular Hopf algebra. Then $H^{(L,J)}$ is a
cosemisimple coalgebra.
\end{theorem}
\begin{remark} This theorem is a generalization
of Corollary \ref{twisted Hopf alg}. Namely, it becomes
Corollary \ref{twisted Hopf alg}  if $L=J$.
\end{remark}

\begin{proof}
Let  $V:= S(Q_L)$ and $W := Q_J^{-1}$. We will show that
\begin{equation}
\label{twisted pairing} \psi(g\otimes h) = \la \lambda,\,h W S(g)
V \ra, \qquad h,g\in H^{(L,J)},
\end{equation}
is a coseparability pairing for $H^{(L,J)}$. To this end we need to
check that $\psi$ satisfies the conditions of
Equation~\eqref{cosep pairing}.

First, we check that $\psi(\tilde{\Delta}(h)) = \eps(h)$ for all $h\in H$ :
\begin{eqnarray*}
\psi(\tilde{\Delta}(h))
&=& \la \lambda, L^{-(1)} h_{\tilde{(1)}} J\I J^{-(1)} S(J^{-(2)})
     S(L^{-(2)} h_{\tilde{(2)}} J\II) S(L\II)S^2(L\I) \ra \\
&=& \eps(h) \la \lambda, L^{-(1)} S(L^{-(2)}) S(L\II)S^2(L\I)\ra  \\
&=& \eps(h) \la \lambda, L\I L^{-(1)} S(L^{-(2)}) S(L\II) \ra \\
&=& \eps(h),
\end{eqnarray*}
by Remark~\ref{nakayama}.

Next, we compute, for all $h,g\in H$ :
\begin{eqnarray*}
\psi (h \otimes g_{\tilde{(1)}}) g_{\tilde{(2)}}
&=& \la \lambda, h W S(g_{\tilde{(1)}}) V\ra  g_{\tilde{(2)}} \\
&=& \la \lambda, h W S( L^{-(1)} g\1 J\I) V\ra   L^{-(2)} g\2  J\II \\
&=& \la \lambda, h W S(J\I) S(g\1) S( L^{-(1)})V \ra   L^{-(2)} g\2  J\II,
\end{eqnarray*}
and also
\begin{eqnarray*}
h_{\tilde{(1)}} \psi(h_{\tilde{(2)}} \otimes g)
&=& L^{-(1)} h\1 J\I  \la \lambda, L^{-(2)} h\2  J\II W S(g) V \ra \\
&=& L^{-(1)} h\1 J\I  \la \lambda,  h\2  J\II W S(g) V S^2(L^{-(2)})\ra  \\
&=& L^{-(1)} S^{-1}[ (J\II W S(g) V S^2(L^{-(2)}))\1] J\I \\
& & \times\,    \la \lambda,  h  (J\II W S(g) V S^2(L^{-(2)}))\2 \ra \\
&=& \la \lambda,  h {J\II}\2 W\2 S(g\1) V\2 S^2({L^{-(2)}}\2) \ra \\
& &  \times\,
  L^{-(1)} S({L^{-(2)}}\1) S^{-1}(V\1) g\2 S^{-1}(W\1) S^{-1}({J\II}\1)J\I.
\end{eqnarray*}
where we used Equation~\eqref{left invariance} and Remark~\ref{nakayama}.

Comparing the results of the above computations and using the non-degeneracy
of $\lambda$, we conclude that the following equations
\begin{equation}
\label{conditions on L and J}
\begin{split}
W S(J\I) \otimes J\II &= {J\II}\2 W\2\otimes S^{-1}(W\1) S^{-1}({J\II}\1)J\I \\
S( L^{-(1)})V \otimes L^{-(2)} &= V\2 S^2({L^{-(2)}}\2) \otimes
                                 L^{-(1)} S({L^{-(2)}}\1) S^{-1}(V\1)
\end{split}
\end{equation}
imply the identity
\begin{equation}
\label{what we need from psi}
\psi (h \otimes g_{\tilde{(1)}}) g_{\tilde{(2)}} =
h_{\tilde{(1)}} \psi(h_{\tilde{(2)}} \otimes g).
\end{equation}

 From Lemma~\ref{twist identities} we see that
equations~\eqref{conditions on L and J} are equivalent to
\begin{equation}
\label{final conditions on L and J}
\begin{split}
\Delta(W) &= J (W \otimes W) (S \otimes S)(J_{21})\\
\Delta(S^{-1}(V^{-1})) &= L^{-1}  (S^{-1}(V^{-1}) \otimes S^{-1}(V^{-1}))
                          (S \otimes S)(L_{21}^{-1}).
\end{split}
\end{equation}
These, in turn, are equivalent to Equation~\eqref{Delta Q}, therefore
the proof is complete.
\end{proof}

\begin{remark} \label{dropped} Suppose that $L=1$. Then it is clear from the
proof of Theorem \ref{HLJ is cosemisimple}
that the unimodularity assumption can be
dropped. In this special case, Theorem
\ref{HLJ is cosemisimple} becomes a known result of Blattner and Montgomery
\cite[Theorem 7.4.2, part 2]{M}, in the case when the algebra
$A$ is one-dimensional.
\end{remark}

\begin{remark} It is clear that $H^{(L,J)}=(H^L)^{(L^{-1}J)}$.
Thus, another proof of Theorem \ref{HLJ is cosemisimple} can be
obtained by combining Corollary \ref{twisted Hopf alg} with
\cite[Theorem 7.4.2, part 2]{M}.
\end{remark}

\begin{remark}
In the case when $k=\mathbb{C}$, the field of complex numbers,
and $H =\mathbb{C}[G]$, the group algebra of a finite group $G$,
the algebraic structure of $(H^J)^*$ was described in \cite[Theorem 3.2]{EG2}.
Below we give a similar description of the algebraic structure of
$(H^{(L,J)})^*$. Obviously, this structure depends only on gauge
equivalence classes of $L, J$.
\end{remark}

Recall \cite{EG1, Mov} that for any twist $J$ of $H =\mathbb{C}[G]$
there exists a subgroup $K\subset G$ and a twist $J^x \in \mathbb{C}[K]
\otimes \mathbb{C}[K]$ which is minimal for $\mathbb{C}[K]$ and is
gauge equivalent to $J$. Such a twist defines a projective representation
$V$ of $K$ with $\dim(V) = |K|^{1/2}$.

Let $L, J$ be twists of this type and $K_L, K_J \subset G$
be the corresponding subgroups for which they are minimal.
We have $L\in \mathbb{C}[K_L] \otimes \mathbb{C}[K_L]$ and
$J\in \mathbb{C}[K_J] \otimes \mathbb{C}[K_J]$.
Let $V_L$ and $V_J$ be the corresponding projective representations.
Let $Z$ be a $(K_L, K_J)$ double  coset in $G$ and let $g\in Z$. Let
$M_g := K_L \cap gK_J g^{-1}$ and define embeddings
$\theta_L : M_g \to K_L$ and  $\theta_J : M_g \to K_J$
by setting $\theta_L(a)=a$ and $\theta_J(a) = g^{-1}ag$ for all $a\in M_g$.
Denote by $W_L$ (respectively, $W_J$) the pullback of the projective
representation $V_L$ (respectively, $V_J$) to $M_g$ by means of
$\theta_L$ (respectively, $\theta_J$).

Let $H_Z^* := \oplus_{g\in Z}\, \mathbb{C}\delta_g \subseteq H^*$, where
delta-functions $\{\delta_g\}_{g\in G}$ form a basis of $H^*$.
Then $H_Z^*$ is a subalgebra of $(H^{(L,J)})^*$ and
$(H^{(L,J)})^* = \oplus_Z\, H_Z^*$ as algebras.

Thus, to find the algebra structure of $(H^{(L,J)})^*$ it suffices
to find the algebra structure of each $H_Z^*$.
The following result is a straightforward generalization of
\cite[Theorem 3.2]{EG2}.

\begin{theorem}
Let $W_L, W_J$ be as above and let $(\widehat{M}_g, \widehat{\pi}_W)$
be any linearization of the projective representation $W := W_L\otimes W_J$
of $M_g$. Let $\zeta$ be the kernel of the projection $\widehat{M}_g\to M_g$,
and let $\chi :\zeta \to \mathbb{C}^\times$ be the character by which
$\zeta$ acts in $W$. Then there exists a one-to-one correspondence
between isomorphism classes of irreducible representations of $H_Z^*$
and isomorphism classes of irreducible representations of $\widehat{M}_g$
with $\zeta$ acting by $\chi$. If a representation $Y$ of $H_Z^*$
corresponds to a representation $X$ of $\widehat{M}_g$ then
\begin{equation}
\label{YX eqn}
\dim(Y) = \frac{\sqrt{|K_L| |K_J|}}{|M_g|} \dim(X).
\end{equation}
\end{theorem}
\begin{proof}
The proof is analogous to \cite[Section 4]{EG2}, where
we refer the reader for details. Note that $|K_L|$ and $|K_J|$
are full squares, because of the minimality of the twists $L$ and $J$,
so the right hand side of Equation~\eqref{YX eqn} is an integer.
\end{proof}

\end{section}

%%%%%%%%%%%%%%%%%%%%%%%%%%%%%%%%%%%%%%%%%%%%%%%%%%%%%%%%%%%%%%%%%%%%%%%%%%%%%%%
\begin{section}
{The $H$-module coalgebra associated to a twist}

Let $H$ be a Hopf algebra.
Recall, that a left {\em $H$-module coalgebra} $(C,\, \Delta_C,\, \eps_C)$
is a coalgebra which is also a left $H$-module via
$h\otimes c \mapsto h \cdot c$, such that
\begin{equation}
\label{Movshes' coalgebra}
\Delta_C(h\cdot c) = \Delta(h) \cdot \Delta_C(c),
\quad \mbox{ and } \quad
\eps_C(h\cdot c) = \eps(h) \eps_C(c),
\end{equation}
for all $h\in H$ and $c\in C$.

Clearly, if $H$ is finite-dimensional, then $C$ is a left
$H$-module coalgebra if and only if $A := C^*$ is a left {\em
$H^*$-comodule algebra} \cite[4.1.2]{M},  i.e., if and only if
the algebra $A$ is a left $H^*$-comodule via $ a \mapsto
\gamma(a)$, such that
\begin{equation*}
\gamma(ab) =\gamma(a)\gamma(b)
\quad \mbox{ and } \quad
\gamma(1_A) = 1 \otimes 1_A
\end{equation*}
for all $a,b \in A$.

It is well known in the theory of Hopf algebras that cleft
$H$-extensions \cite[7.2]{M} of algebras are precisely cocycle
crossed products with $H$ \cite[Theorem 7.2.2]{M}. Since a twist for a
finite-dimensional $H$ is naturally a $2-$cocycle for $H^*$, it is
quite clear that twists for $H$ can be characterized in terms of
extensions.

We will need the following definition of a Galois extension, cf.\
\cite[8.1.1]{M}.

\begin{definition}
\label{galois}
Let $A$ be a left $H$-comodule algebra with the structure map
\linebreak $\gamma : A \to H \otimes A$. Then the $H$-extension of algebras
$A^{\text{co}\,H}\subset A$
is left {\em $H$-Galois} if the map $A \otimes_{A^{\text{co}\,H}} A
\to H \otimes_k A$ given by $a\otimes b \mapsto (1\otimes a)\gamma(b)$
is bijective.
\end{definition}

Let $H$ be a finite-dimensional Hopf algebra and $J\in H\otimes H$
be a twist for $H$. Then the left regular $H$-module $C:=H^{(J)}=
H$ with comultiplication and counit
\begin{equation}
\Delta_{H^{(J)}}(c) = \Delta(c) J, \qquad
\eps_{H^{(J)}}(c) =\eps(c), \qquad c\in C,
\end{equation}
is a left $H$-module coalgebra. Moreover, the corresponding
$H^*$-extension $k\subset (H^{(J)})^*$ is Galois, since $J$ is invertible.

Let $J$ and $J^x = \Delta(x) J (x^{-1} \otimes x^{-1})$ be gauge
equivalent twists, where $x\in H$ is invertible and is such that $\eps(x)=1$.
Then $H^{(J^x)}=H$ with the comultiplication $\Delta_{H^{(J^x)}}(c)
= \Delta(c)J^x$ is another left $H$-module coalgebra such that the map
$c \mapsto cx$ is an $H$-module coalgebra isomorphism between $H^{(J)}$
and $H^{(J^x)}$.

Conversely, let $C$ be an $H$-module coalgebra and let $i : C \to
H$ be an isomorphism of left $H$-modules (where $H$ is viewed as a
left regular module over itself), such that the corresponding
$H^*$-extension $k\subset C^*$ is Galois. Then one can check by a
direct computation that
\begin{equation}
\label{formula for J}
J := (i \otimes i)\Delta_C(i^{-1}(1)) \in
H\otimes H
\end{equation}
is a twist for $H$ and $i : C \cong H^{(J)}$ is an
$H$-module coalgebra isomorphism. Here the Galois property above
is equivalent to $J$ being invertible.

\begin{proposition}\label{1to1}
The above two assignments define a one-to-one correspondence
between:
\begin{enumerate}
\item[1.]
gauge equivalence classes of twists for $H$, and
\item[2.]
isomorphism classes of $H$-module coalgebras $C$ isomorphic to the
regular $H$-module and such that the corresponding
$H^*$-extension $k\subset C^*$ is Galois.
\end{enumerate}
\end{proposition}

\begin{corollary}
\label{semisimplicity all around} If $H$ is a finite-dimensional
cosemisimple Hopf algebra, then any $H$-module
coalgebra $C$ isomorphic to the
regular $H$-module and such that the corresponding
$H^*$-extension $k\subset C^*$ is Galois, is cosemisimple.
%as in Equation~\eqref{Movshes' coalgebra}
\end{corollary}

\begin{proof} This follows from Theorem~\ref{HLJ is cosemisimple}
and Remark \ref{dropped}.
\end{proof}

\end{section}

%%%%%%%%%%%%%%%%%%%%%%%%%%%%%%%%%%%%%%%%%%%%%%%%%%%%%%%%%%%%%%%%%%%%%%%%%%%%%
\begin{section}{Non-degenerate twists and non-commutative
``lowering of indices''}

Let $H$ be a finite-dimensional Hopf algebra over $k$.

\begin{definition}
\label{nondeg twist}
We will say that a twist $J$ for $H$ is {\em non-degenerate} if
the corresponding coalgebra $H^{(J)}$ is simple.
\end{definition}

\begin{remark}
If $H$ is a group algebra over a field of
characteristic prime to the order of the group
then the notion of a non-degenerate
twist coincides with the notion of a minimal twist
from \cite{EG1} (i.e., this property is equivalent to the triangular Hopf
algebra $ (H^J,J_{21}^{-1}J)$ being minimal).
\end{remark}

\begin{remark}
\label{Sweedler}
A Hopf algebra possessing
a non-degenerate twist need not be semisimple or cosemisimple. For
example, in \cite{AEG} it was explained that Sweedler's
4-dimensional Hopf algebra $H$ (which is neither semisimple nor
cosemisimple) possesses a 1-parameter family of twists
$J(t):=1\otimes 1-(t/2)gx\otimes x$ (where $g$ is the group-like element
and $x$ the skew primitive element). It is easy to show that
$J(t)$ is non-degenerate for $t\ne 0$. Indeed, it is easy to prove
by a direct calculation that the coalgebra $H^{(J(t))}$ does not
have group-like elements. But any 4-dimensional coalgebra without
group-like elements is necessarily simple.
\end{remark}

Now consider a twist $J\in H\otimes H$. Observe that $H^*_{(J)}:=
(H^{(J)})^*$ is a right $H$-module algebra in a natural way:
\begin{equation}
\la a \cdot h,\, c \ra = \la a,\, hc\ra, \qquad a\in
H^*_{(J)},\,c\in H^{(J)},\,h\in H.
\end{equation}
The algebra  $H^*_{(J)}$ is simple if and only if $J$ is
non-degenerate. In this case the Skolem-Noether theorem for Hopf
algebras \cite{Mas} (see also \cite[6.2.4]{M}) says that the
action of $H$ is inner, i.e., there is a convolution invertible
map ${\pi}\in \Hom_k(H,H^*_{(J)})$ such that
\begin{equation}
\label{inner action}
a \cdot h = \overline{{\pi}}(h\1) a {\pi}(h\2),
\qquad \mbox{for all } a\in H^*_{(J)},\, h\in H,
\end{equation}
where $\overline{{\pi}}$ is the convolution inverse for ${\pi}$.
Furthermore, we may (and will) assume that ${\pi}(1)=
\overline{{\pi}}(1) =1$. Such a map $\pi$ will be called a {\em
Skolem-Noether map}.

It is easy to check that a Skolem-Noether map is unique up to a
gauge transformation, i.e., that two Skolem-Noether maps $\pi,\pi'$ are
linked by the relation $\pi'(h)=\eta(h\1)\pi(h\2)$, where $\eta\in
H^*$ is invertible and $\eta(1)=1$.

\begin{remark} If $H$ is the group algebra of a group $G$ then
a Skolem-Noether map is the same thing as an irreducible projective
representation of $G$. Thus, the theory of Skolem-Noether maps
is a generalization of the theory of projective representations
to Hopf algebras.
\end{remark}

Comparing the values of $(a\cdot g)\cdot h$ and $a\cdot (gh)$, it
is easy to deduce from \linebreak Equation~\eqref{inner action}
that $c (h\otimes g) := \pi(h\1 g\1) \overline{\pi}(g\2) \overline{\pi}(h\2)$
is a scalar for all $g,h \in H$. Furthermore, the form $c\in
H^*\otimes H^*$ is an invertible element (by convolution
invertibility of $\pi$), and is a $2$-cocycle for $H$, i.e.,
\begin{equation}
%\label{cocycle eqn}
c(f\1g\1\otimes h) c(f\2\otimes g\2) =
c(f \otimes g\1h\1) c(g\2\otimes h\2), \qquad \mbox{for all } f,g,h \in  H.
\end{equation}
Indeed, we have
\begin{eqnarray*}
\lefteqn{ c(f\1g\1\otimes h) c(f\2\otimes g\2) =} \\
&=& \pi(f\1g\1h\1) \overline{\pi}(h\2) \overline{\pi}(f\2g\2)
    \pi(f\3g\3)  \overline{\pi}(g\4) \overline{\pi}(f\4) \\
&=& \pi(f\1g\1h\1)  \overline{\pi}(h\2) \overline{\pi}(g\2)
\overline{\pi}(f\2),\\ \lefteqn{c(f \otimes g\1h\1) c(g\2\otimes
h\2) =}\\ &=& \pi(f\1g\1h\1) \overline{\pi}(g\2h\2) c(g\3\otimes
h\3) \overline{\pi}(f\2) \\ &=& \pi(f\1g\1h\1)
\overline{\pi}(g\2h\2) {\pi}(g\3h\3) \overline{\pi}(h\4)
    \overline{\pi}(g\4)  \overline{\pi}(f\2)\\
&=& \pi(f\1g\1h\1)  \overline{\pi}(h\2) \overline{\pi}(g\2) \overline{\pi}(f\2).
\end{eqnarray*}

Thus, $c$ defines a twist for $H^*$.
We will denote this twist by $D_\pi(J)$ (the dual twist to $J$,
constructed
via $\pi$).

It is easy to check that while the twist $D_\pi(J)$ does depend
on the choice of $\pi$, its gauge equivalence class does not.
Namely, for a fixed $J$, the map $\pi\to D_\pi(J)$ commutes with
gauge transformations in $H^*$.

Let $NT(H)$, $NT(H^*)$ be the sets of gauge equivalence classes of
non-degenerate twists for $H$, $H^*$, respectively. Below we
explain that the map $D_H: J \to D_\pi(J)$ establishes a duality
between $NT(H)$ and $NT(H^*)$.

\begin{remark} Of course, both $NT(H)$ and $NT(H^*)$ may well be
empty. It is clearly so, for example, if the dimension of $H$ is
not a square.
\end{remark}

We adopt the viewpoint of Miyashita-Ulbrich actions \cite{DT2, U},
which we recall next.  Let $H$ be a finite-dimensional Hopf
algebra and $R$ be a simple algebra. Suppose there is a right
$H$-Galois coaction on $R$, $y\mapsto y\0\otimes y\1$, by which we
mean an $H$-coaction making $R$ right $H$-Galois over $k$ (i.e.
$R$ is an $H-$comodule algebra and $R^{coH}=k$). Then we have the
right $H$-action (the Miyashita-Ulbrich action) $x\ract h,\, x\in
R,\, h\in H$ characterized by
\begin{equation}
\label{MU action}
xy = y\0 (x \ract y\1), \qquad x,y\in R,
\end{equation}
see \cite[Theorem 3.4(i)]{DT2}. The construction of the left
Miyashita-Ulbrich action from a left $H$-Galois coaction on $R$
is completely similar.

\begin{remark}
The Miyashita-Ulbrich actions can be defined in a more general situation,
when $k$ is  an arbitrary (not necessarily algebraically closed) field
and $R$ is an Azumaya (or finite-dimensional central simple) algebra \cite{DT2}.
\end{remark}

The map $D_H$ can be clearly explained in terms of the
Miyashita-Ulbrich actions as follows.
%Let again $R$ be a simple
%algebra with a right Galois coaction  $y\mapsto y\0\otimes y\1\,
%(y\in R)$ of a finite-dimensional Hopf algebra $H$ and
Let $\ract$ be the corresponding right Miyashita-Ulbrich action of
$H$ \eqref{MU action}. The left $H^*$-coaction $x \mapsto
x_{(-1)}\otimes x_{(0)}\, (x\in R)$ adjoint to this action is a
left $H^*$-Galois coaction on $R$ by \cite[Satz 2.8]{U},
\cite[Theorem 5.1]{DT2}. The left Miyashita-Ulbrich action $t\lact
x\, (t\in H^*, x\in R)$ arising from the last $H^*$-coaction
coincide with the one that is adjoint to the original
$H$-coaction, since we have
\begin{equation*}
y\0 \la x_{(-1)},\, y\1 \ra x\0 = y\0 (x\ract y\1) = xy = (x_{(-1)}\lact y)x\0,
\end{equation*}
and hence, $y\0 \la t,\, y\1 \ra = t \lact y$ for all $t\in H^*$ and $y\in R$,
since $x\mapsto x_{(-1)} \otimes x\0$ is Galois,
see \cite[(3.1)]{DT2}. Thus we have a natural 1-1 correspondence between the right
$H$-Galois coactions on $R$ and the left $H^*$-Galois coactions on $R$.
This can be regarded as a very special case of \cite[Theorem 6.20]{DT2}.
As is easily seen, this induces a 1-1 correspondence between the sets
of respective equivalence classes of coactions, where two right coactions
are equivalent if the corresponding comodule algebras are isomorphic.
By Proposition~\ref{1to1}, this gives a bijective map between
$NT(H)$ and $NT(H^*)$.

\begin{theorem}
\label{bijection for twists}
The map $D_H: NT(H)\to NT(H^*)$ is a bijection and $D_H^{-1}=D_{H^*}$.
\end{theorem}
\begin{proof}
Write $R= H^*_{(J)}$. The map $\pi : H\to R$ is left $H^*$-linear,
where  $H^*$ acts on $R$ by the Miyashita-Ulbrich action, since we have
\begin{eqnarray*}
\lefteqn{ (x_{(-1)} \lact \pi(h)) x\0 = x\pi(h) =\pi(h\1)(x \cdot h\2) }\\
&=& \pi(h\1) \la x_{(-1)},\,h\2\ra x\0
= \pi(x_{(-1)} \actl h) x\0,\qquad h\in H,\,x\in R.
\end{eqnarray*}
It follows from \cite[Theorems 9,11]{DT1} that $R$ is cleft right $H$-Galois
with coinvariants $k$ (by counting dimension), and
\begin{equation}
\pi : H_{(D_\pi(J))} \to R
\end{equation}
is an isomorphism of right $H-$comodule algebras.
 From formula \eqref{formula for J} we see that the bijective map between
$NT(H)$ and $NT(H^*)$ mentioned above coincides with $D_H$. It
also follows from the construction of this map that
$D_H^{-1}=D_{H^*}$.
\end{proof}

\begin{remark}
When $H$ is semisimple and cosemisimple,
the set $NT(H)$ is finite. This follows from the result
of H.-J.~Schneider \cite{Sch} stating that the number
of isomorphism classes of semisimple $H$-comodule algebras
of fixed dimension is finite.
\end{remark}

\begin{example}
Let $F$ be a finite field of characteristic $p\neq 2$,
and $\psi: F\to \mathbb{C}^\times$ a nontrivial additive character
(for instance, $F=\mathbb{Z}/p\mathbb{Z}$, $\psi(n)=e^{2\pi in/p}$).
Let $V$ be a finite-dimensional vector space over $F$.
Suppose that $H=\mathbb{C}[V]$.
In this case $NT(H)$ can be identified with the
set of non-degenerate alternating bilinear forms $\omega$ on
$V^*$, by the assignment $\omega\mapsto J_\omega$, where
$J_\omega$ is defined by
\begin{equation}
((\psi\circ \chi_1)\otimes (\psi\circ\chi_2))(J_\omega)=
\psi(\frac{1}{2}\omega(\chi_1,\chi_2)),\qquad \chi_i\in V^*.
\end{equation}
On the other hand, $H^*$ is isomorphic to $\mathbb{C}[V^*]$,
via $\chi \to \psi\circ \chi\in H^*$.
Therefore, $NT(H^*)$ is identified with the set of non-degenerate alternating
bilinear forms on $V$. Thus, it is easy to guess the canonical
map $D_H: NT(H)\to NT(H^*)$ in this case. Namely, it is not hard
to check that this map is the usual inversion of a non-degenerate form,
which in  classical tensor calculus is often called the ``lowering
of indices'' \cite{DFN}. This motivated the title of this subsection.
\end{example}

\begin{example}
This example is a generalization of the previous one.
Now assume that $H=k[G]$, where $G$ is a finite group
of order coprime to $\text{char}(k)$. In this case,
$NT(H^*)$ is the set of cohomology classes of 2-cocycles
on $G$ with values in $\mathbb{C}^\times$,
which are non-degenerate in the sense
of \cite{Mov}. Furthermore, the bijection $D_H$ in this case is
nothing but the bijection from \cite{Mov} between
such cohomology classes and classes of non-degenerate twists for
$k[G]$,
which plays a central role in the group theoretic description of
twists in $k[G]$. Thus, the theory of this subsection is a
generalization of the theory of \cite{Mov} and \cite{EG1}
to the non-cocommutative case.
\end{example}

\begin{example}
Let us give an example of a non-degenerate twist for a Hopf algebra which
does not come from group theory. Let $A$ be a
finite-dimensional Hopf algebra and
$\{a_i\}$ and $\{\phi_i\}$ be dual bases in $A$ and $A^*$ respectively.
Then it is straightforward to check that
\begin{equation}
\label{J in the double}
J := \sum_i\, (\phi_i\otimes 1) \otimes (\eps \otimes a_i)
\end{equation}
is a twist for the Hopf algebra $H := A^{*\text{op}} \otimes A$.
Then the twisted Hopf algebra is $H^J\cong D(A)^{*op}$, the
(opposite) dual of
the Drinfeld double of $A$. Furthermore, it was proved in
\cite[Theorem 6.1]{Lu} that
$(H^{(J)})^* \cong A^*\# A$, where $A$ acts on its dual via $\actl$.
This algebra $A^*\# A$ is called the Heisenberg double in the literature,
and is known to be simple \cite[Corollary 9.4.3]{M}. Thus, $J$
is a non-degenerate twist for $H$.
\end{example}

\begin{remark}
Hopf algebras with a non-empty set of non-degenerate twists
may be called ``quantum groups of central type'' by analogy with
\cite{EG1}. It would be interesting to find more examples of
such Hopf algebras.
\end{remark}

\end{section}

%%%%%%%%%%%%%%%%%%%%%%%%%%%%%%%%%%%%%%%%%%%%%%%%%%%%%%%%%%%%%%%%%%%%%%%%%%%%%%%
\begin{section} {The classification of finite-dimensional
triangular cosemisimple Hopf
algebras in positive characteristic}

The classification of isomorphism classes of finite-dimensional
triangular cosemisimple Hopf algebras over $\mathbb{C}$ was
obtained in \cite{EG1} (since by \cite{LR}, such Hopf algebras are
also semisimple). Namely, in \cite{EG1} they were put in bijection
with certain quadruples of group-theoretical data. Our next goal
is to classify finite-dimensional triangular cosemisimple Hopf
algebras over a field $k$ of characteristic
$p>0$ in similar terms.

We start by defining the class of quadruples of group-theoretical
data which is suitable in the positive characteristic case.

\begin{definition}\label{trquad}
Let $k$ be a field of characteristic $p>0.$ A quadruple
$(G,K,V,u)$ is called a {\em triangular quadruple (over $k$)} if $G$ is a
finite group, $K$ is a $p'-$subgroup of $G$
(i.e., a subgroup of order coprime to $p$), $V$ is an irreducible
projective representation of $K$ of dimension $\sqrt{|K|}$ over $k,$
and $u\in G$ is a central element of order
less or equal than $\min(2,p-1)$ (so if $p=2$ then $u=1$).
\end{definition}

The notion of an isomorphism between triangular quadruples is clear.
Given a triangular quadruple $(G,K,V,u),$ one can construct a
triangular cosemisimple Hopf algebra $H(G,K,V,u)$ in the following
way. As in \cite{EG1}, first construct the twist $J(V)$ of $k[K]$
corresponding to $V$ (well defined up to gauge transformations),
and then define the triangular Hopf algebra
$$H(G,K,V,u):=(k[G]^{J(V)},J(V)_{21}^{-1}J(V)R_u),$$ where
$R_u=1$ if $u=1$ and
$$R_u:=\frac{1}{2}(1\ot 1+1\ot u+u\ot 1-u\ot u)$$
if $u$ is of order $2$ (note that the $1/2$ makes sense, since by the
definition, if $u$ has order $2$ then $p$ is odd).
As a Hopf algebra, $H(G,K,V,u)$ is isomorphic to $k[G]^{J(V)}$, and
therefore it is cosemisimple by Corollary \ref{twisted Hopf alg}
(as $k[G]$ is obviously unimodular).

Conversely, let $(H,R)$ be a finite-dimensional triangular
cosemisimple Hopf algebra over $k.$ Let $(H_{\text{min}},R)$ be
the minimal triangular sub Hopf algebra of $(H,R).$ Since
$H_{\text{min}}\cong H_{\text{min}}^{*\text{op}}$ as Hopf algebras,
$H_{\text{min}}$ is cosemisimple and semisimple. Hence by
\cite{EG1}, there exists a finite $p'-$group $K$ and a twist $J$
for $K$ such that $$ (H_{\text{min}},RR_u)\cong
(k[K]^J,J_{21}^{-1}J), $$ where $u\in H_{\text{min}}$ is a central
group-like element satisfying $u^2=1,$ and $R_u$ is as above. Note
that $u$ is the Drinfeld element of $H_{\text{min}},$ hence of $H$
as well.

Consider now the triangular Hopf algebra
$(H,R)^{J^{-1}}=(H^{J^{-1}},R^J),$ and set $B:=H^{J^{-1}}.$ By
\cite[Theorem 1.3.6]{G1}, $H$ is unimodular and hence by
Corollary \ref{twisted Hopf alg}, $B$ is cosemisimple.

Our next goal is to show that $u$ is central in $H$. If $u=1$,
there is nothing to prove. Thus, we may assume that $u$ has order
$2$, which in particular implies that $p$ is odd
(since the order of $u$ divides the dimension of $H_{\text{min}}$, while
$p$ does not divide this dimension, as $H_{\text{min}}$ is semisimple
and cosemisimple).

Let $\mathcal{B}$ be the cocommutative Hopf superalgebra
corresponding to $(B,u)$ by \cite[Corollary 3.3.3]{AEG} (this
corollary was proved for characteristic zero, but applies verbatim
in odd characteristic). In particular $u\in \mathcal{B}$ acts by
parity.

\begin{lemma}\label{cosssup}
The Hopf superalgebra $\mathcal{B}$ is cosemisimple.
\end{lemma}

\begin{proof}
Form the biproduct $\overline{\mathcal{B}}:=k[\mathbb{Z}_2]\ltimes
\mathcal{B}$ where the generator $g$ of $\mathbb{Z}_2$ acts on
$\mathcal{B}$ by parity; it is an ordinary Hopf algebra (see
\cite{AEG}). Consider the algebra $\overline{\mathcal{B}}^*.$ On
one hand, $\overline{\mathcal{B}}^*=k[\mathbb{Z}_2]\ltimes
\mathcal{B}^*$ as an algebra, and on the other hand,
$\overline{\mathcal{B}}^*=k[\mathbb{Z}_2]\ot B^*$ as an algebra.
Since $B$ is
cosemisimple, the result follows.
\end{proof}

Since $\mathcal{B}$ is cocommutative and cosemisimple, it is a
group algebra $k[G]$ of some finite group $G.$ In particular,
$\mathcal{B}$ is purely even, and hence $u$ is central in $B$ and
hence central in $H$ as well.

We conclude that $(H,R)^{J^{-1}}=(H^{J^{-1}},R_u)=(k[G],R_u).$ But
since $k[K]$ is embedded in $k[G]$ as a Hopf algebra, $K$ is a
subgroup of $G.$ We thus assigned to $(H,R)$ a triangular
quadruple $(G,K,V,u).$

In fact we have the following theorem.

\begin{theorem} \label{poschar}
Let $k$ be a field of characteristic $p>0.$
The above two assignments define a one to one correspondence
between:

\begin{enumerate}

\item isomorphism classes of finite-dimensional
triangular cosemisimple Hopf algebras
 over $k,$ and

\item isomorphism classes of triangular quadruples $(G,K,V,u)$
over $k$.

\end{enumerate}
\end{theorem}

\begin{proof} It remains to show that the two assignments are inverse to
each other. Indeed, this follows from the results of \cite[Section
5]{EG1}.
\end{proof}

\begin{remark} This theorem and Remark \ref{involutive} imply that in a
finite-dimensional cosemisimple triangular Hopf algebra, one has
$S^2=\id$. So Kaplansky's conjecture that in a semisimple Hopf
algebra, one has $S^2=\id$, is valid in the triangular case.
\end{remark}

As a corollary we are now able to classify twists for finite
groups in positive characteristic.

\begin{theorem} \label{poschartw} Let $k$ be a
field of characteristic $p>0.$ Let $G$ be a finite group. There is
a one to one correspondence between:

\begin{enumerate}
\item gauge equivalence classes of twists for $k[G],$ and

\item pairs $(K,V)$ where $K$ is a $p'-$subgroup of $G$ and $V$ is
an irreducible projective representation of $K$ of dimension
$\sqrt{|K|}$ over $k,$ modulo inner automorphisms of $G$.

\end{enumerate}
\end{theorem}
\begin{proof}
Follows from Theorem \ref{poschar} since for any twist $J$ for
$k[G],$ $(k[G]^J,J_{21}^{-1}J)$ is cosemisimple triangular by
Corollary \ref{twisted Hopf alg}.
\end{proof}

\begin{remark}
Theorem~\ref{poschartw} implies that the representation theory of $(k[G]^J)^*$
for twists $J \in k[G]\otimes k[G]$ is described exactly as in \cite{EG2}.
\end{remark}

\end{section}
%%%%%%%%%%%%%%%%%%%%%%%%%%%%%%%%%%%%%%%%%%%%%%%%%%%%%%%%%%%%%%%%%%%%%%%%%%%%%%
%%%%%%%%%%%%%%%   BIBLIOGRAPHY %%%%%%%%%%%%%%%%%%%%%%%%%%%%%%%%%%%%%%%%%%%%%%%
%%%%%%%%%%%%%%%%%%%%%%%%%%%%%%%%%%%%%%%%%%%%%%%%%%%%%%%%%%%%%%%%%%%%%%%%%%%%%

\bibliographystyle{ams-alpha}

\end{document}